\documentclass[11pt]{article}
\usepackage{graphicx}
\usepackage{color}
\usepackage{multirow}
\usepackage{lscape}
\usepackage{amssymb,latexsym,amsmath}
\usepackage{hyperref}
\usepackage{epstopdf}
\newtheorem{lemma}{Lemma}[section]
\newtheorem{theo}[lemma]{Theorem}

\newtheorem{res}[lemma]{Result}

\newcommand{\proof}{\noindent{\em Proof: }}

\newcommand{\forme}[1]{}

\def\wbull{\hfill\vrule height .9ex width .8ex depth -.1ex}

\DeclareMathOperator{\PG}{PG}
\def\St{{\mathsf{Star}}}
\def\Li{{\mathsf{Line}}}

\begin{document}

\date{\today}
\title{Cameron-Liebler line classes in $\PG(3,5)$}
\author{{\bf Alexander L. Gavrilyuk}\\
Center for Math Research and Education,
Pusan National University,\\
2, Busandaehak-ro 63beon-gil, Geumjeong-gu, Busan, 46241, Republic of Korea\\
and\\
Krasovskii Institute of Mathematics and Mechanics,\\ 
Kovalevskaya str., 16, Ekaterinburg 620990, Russia\\
e-mail: alexander.gavriliouk@gmail.com\\
\\
{\bf Ilia Matkin}\\
Faculty of Mathematics, Chelyabinsk State University,\\
Kashirinykh str., 129, Chelyabinsk 454001, Russia\\
e-mail: ilya.matkin@gmail.com}
\maketitle

\begin{abstract}
We complete a classification of Cameron-Liebler line classes in $\PG(3,5)$, 
and show in a uniform way all non-existence results for those in $\PG(3,q)$, $q\leq 5$.
\end{abstract}

\section{Introduction}

Let $\PG(n,q)$ denote the $n$-dimensional projective space over the finite field $\mathbb{F}_q$ 
with $q$ elements.
Let $A$ be the point-line incidence matrix of $\PG(n,q)$, i.e., the rows of $A$ are indexed by 
the set of points of $\PG(n,q)$, its columns are indexed by the set of lines, and $(A)_{p,\ell}=1$ 
if $p\in \ell$, and $(A)_{p,\ell}=0$ otherwise.
We consider $A$ as a matrix over $\mathbb{Q}$, and let ${\sf R}(A)$ denote the row space of $A$.

A set $\mathcal{L}$ of lines of $\PG(n,q)$ is called a {\it Cameron-Liebler line class} 
\cite{DrudgeThesis}, \cite{Penttila}
if its characteristic vector $\chi_{\mathcal{L}}$ satisfies $\chi_{\mathcal{L}}\in {\sf R}(A)$. 
In their study \cite{CameronLiebler}
of collineation groups of $\PG(n,q)$, $n\geq 3$, that have equally many orbits 
on lines and on points, Cameron and Liebler showed that the characteristic vector 
of a line orbit of such a group should enjoy this property. 

One can see that an empty set of lines, the set $\Li(H)$ of all lines in a hyperplane $H$ or, dually, 
the set $\St(P)$ of all lines through a point $P$, and $\St(P)\cup\Li(H)$ provided that $P\notin H$ 
are examples of Cameron-Liebler line classes. 
We call these examples {\it trivial}.
It was conjectured in \cite{CameronLiebler} (cf. \cite[Section~6]{DrudgeThesis}) 
that, up to their complements, these are the only Cameron-Liebler line classes. 

The conjecture 
turned out to be wrong in $\PG(3,q)$: the first counter example was constructed 
by Drudge \cite{Drudge}, and many more non-trivial line classes have been found during 
the last decade \cite{BruenDrudge,Cossidente,CP,Metsch-and-Co,GMP,GavrilyukMetsch,GovaertsPenttila,FMX,Rodgers}. 
However, its validity in $\PG(n,q)$ with $n>3$ remains an open question.

Cameron-Liebler line classes enjoy a number of properties, which can be taken as their 
equivalent definitions (see \cite{Penttila,DrudgeThesis}). 
In this paper we will focus on $\PG(3,q)$ in which case 
a Cameron-Liebler line class $\mathcal{L}$ (with {\it parameter} $x$) can be defined 
as follows: there exists a non-negative integer $x$ such that every spread $S$ of $\PG(3,q)$ 
shares precisely $x$ lines with $\mathcal{L}$, i.e., $|S\cap \mathcal{L}|=x$. 
It immediately follows from this definition that 
the complementary set $\overline{\mathcal{L}}$ of lines 
is a Cameron-Liebler line class with parameter $q^2+1-x$, so we may further assume 
$x\leq \lfloor\frac{q^2+1}{2}\rfloor$.
For the trivial examples listed above, their parameter $x$ is: 
$0$ for the empty set of lines, $1$ for $\Li(H)$ or $\St(P)$, and 
$2$ for $\Li(H)\cup \St(P)$, and these are the only line classes with parameter $x\leq 2$.

Despite the fact that the conjecture is wrong in $\PG(3,q)$, 
the study of Cameron-Liebler line classes in $\PG(3,q)$ is still of great interest.
We observe that all known non-trivial examples have relatively large parameter $x\simeq q^2/2$, 
although 
the best known lower bound for the parameter $x$ of a non-trivial Cameron-Liebler line 
class is $x\gtrsim q^{4/3}$, see \cite{Metsch2}.
Secondly, for given $q$, one can try to verify the conjecture in $\PG(n,q)$ for all $n>3$ 
provided that a complete list of Cameron-Liebler line classes in $\PG(3,q)$ is known, 
see \cite[Section~6.2]{DrudgeThesis}, \cite{Gavrilyuk} 
(see also \cite{Filmus,CLk} for a higher dimensional generalization of 
Cameron-Liebler line classes).

The complete lists of Cameron-Liebler line classes in $\PG(3,3)$ and in $\PG(3,4)$ 
were obtained in \cite{DrudgeThesis} and \cite{Gavrilyuk}, respectively, see 
also Section \ref{section-comp}.
Up to their complements and a point-plane duality, 
the following non-trivial Cameron-Liebler line classes with parameter $x$ 
are known to exist in $\PG(3,5)$:
\begin{itemize}
\item $\mathcal{G}$ with $x=10$ constructed in \cite{GavrilyukMetsch}; 
\item $\mathcal{R}$ with $x=12$ constructed in \cite{Rodgers} 
(later generalised in \cite{Metsch-and-Co}, \cite{FMX} to an infinite family 
with $x=\frac{q^2-1}{2}$ for $q\equiv 5{\rm~or~}9~({\rm mod~}12)$);
\item $\mathcal{R}^+$ with $x=13$ (there exists a plane $\pi_0$ that is disjoint with $\mathcal{R}$, 
see \cite{RodgersPhD}, and $\mathcal{R}^+=\mathcal{R}\cup \Li(\pi_0)$);
\item $\mathcal{D}$ with $x=13$ (a member of an infinite family with $x=\frac{q^2+1}{2}$ 
constructed in \cite{BruenDrudge});
\item $\mathcal{P}$ with $x=13$ (a member of another infinite family with $x=\frac{q^2+1}{2}$ 
independently constructed in \cite{Cossidente} and in \cite{GMP}).
\end{itemize}

Our main result shows that this list is complete.

\begin{theo}\label{theo-main}
Up to a polarity of $\PG(3,5)$ or taking the complementary set of lines, 
a non-trivial Cameron-Liebler line class in $\PG(3,5)$ 
is projectively equivalent to one of the following: $\mathcal{G}$, $\mathcal{R}$, 
$\mathcal{R}^+$, $\mathcal{D}$, $\mathcal{P}$.
\end{theo} 

This theorem can be used to show that all Cameron-Liebler line classes in $\PG(n,5)$, $n>3$, 
are trivial, see \cite{PGn5}.
We were informed by John Bamberg (a private communication) that 
he found all non-trivial Cameron-Liebler line classes in $\PG(3,5)$ 
with a stabiliser divisible by an element of order $3$, $5$, $13$, or $31$.
His result coincides with our list.

Given a Cameron-Liebler line class $\mathcal{L}$ with parameter $x$ in $\PG(3,q)$, one can construct 
three more line classes associated with $\mathcal{L}$: the complementary set of lines 
$\overline{\mathcal{L}}$ with parameter $q^2+1-x$, the image $\mathcal{L}^{\rho}$ of $\mathcal{L}$ 
under a point-plane duality $\rho$ of $\PG(3,q)$ (a {\it dual} line class) with parameter $x$, and 
a dual line class to $\overline{\mathcal{L}}$ with parameter $q^2+1-x$.
In what follows, we say that $\mathcal{L}$ is {\it unique} if any Cameron-Liebler line 
class $\mathcal{L}'$ with parameter $x$ in $\PG(3,q)$ is projectively equivalent 
to (at least) one of the line classes $\mathcal{L}$, $\overline{\mathcal{L}}$, $\mathcal{L}^{\rho}$, $\overline{\mathcal{L}^{\rho}}$.

The paper is organised as follows. In Section \ref{section-basic}, we recall 
the notion of patterns of lines with respect to a Cameron-Liebler line class, 
which was introduced in \cite{Gavrilyuk}.
Roughly speaking, given a Cameron-Liebler line class $\mathcal{L}$ and a line $\ell$ of $\PG(3,q)$, 
the pattern of $\ell$ w.r.t. $\mathcal{L}$ shows how the lines of $\mathcal{L}$ interface 
with the set of plane pencils containing $\ell$.
The set of all possible patterns may give some insight 
about the structure of a putative line class $\mathcal{L}$, 
sometimes enough to show its existence 
and uniqueness, or its non-existence, see \cite{Metsch-and-Co,GMP,Gavrilyuk,GavrilyukMetsch}. 
Unfortunately, the number of admissible patterns quickly grows as $q$ increases, 
which makes their further ``ad-hoc'' combinatorial analysis very complicated.
In Section \ref{section-main}, we show that such an analysis can be reduced to solving 
a system of Diophantine equations, and demonstrate how it works for $q\in \{3,4,5\}$ 
(thereby we prove in a uniform way all previously known non-existence results 
for Cameron-Liebler line classes in $\PG(3,q)$ with $q\leq 5$).
In Section \ref{section-5}, we apply this approach to determine all non-trivial 
Cameron-Liebler line classes in $\PG(3,5)$, which shows our main result.
In Section \ref{section-comments}, we report some computational results for $q\geq 7$ 
and discuss further problems.


\section{Patterns and the modular equality}\label{section-basic}

We first recall some basic properties of Cameron-Liebler line classes.

\begin{res} [\cite{CameronLiebler,Penttila}]\label{res-basic}
Let $\mathcal{L}$ be a Cameron-Liebler line class with parameter $x$ in $\PG(3,q)$.
\renewcommand{\labelenumi}{\rm (\alph{enumi})}
\begin{enumerate}
\item $|\mathcal{L}|=x(q^2+q+1)$.
\item For every line $\ell$ of $\PG(3,q)$, the number of lines of $\mathcal{L}$ incident to $\ell$ 
equals $x(q+1)+(q^2-1)\chi_{\mathcal{L}}(\ell)$.
\item For any pair of skew lines $\ell,\ell'$ of $\PG(3,q)$, 
the number of lines of $\mathcal{L}$ incident to both $\ell$ and $\ell'$ equals 
$x+q(\chi_{\mathcal{L}}(\ell)+\chi_{\mathcal{L}}(\ell'))$.
\end{enumerate}
\end{res}

Suppose that $\PG(3,q)$ contains a Cameron-Liebler line class $\mathcal{L}$ with parameter $x$.
Pick a line $\ell$ of $\PG(3,q)$. Let $\pi_1$, $\dots$, $\pi_{q+1}$ be the $q+1$ planes 
containing $\ell$, and $P_1$, $\dots$, $P_{q+1}$ be the $q+1$ points on $\ell$.
We define a square matrix $T(\ell)=(t_{ij})$ of size $q+1$ with integer entries given by
\[
t_{ij}:=|((\Li(P_i)\cap\St(\pi_j))\setminus\{\ell\})\cap \mathcal{L}|,~~1\le i,j\le q+1. 
\]

The set consisting of the matrix $T$, and every matrix obtained from this one 
by a permutation of the rows and a permutation of the columns is called the {\it pattern} of $\ell$ 
with respect to ${\cal L}$. (Note that the transpose of $T$ corresponds to the pattern of a line 
with respect to the image of $\mathcal{L}$ under a point-plane duality of $\PG(3,q)$.)
We represent a pattern by one of its matrices. By slight abuse of notation 
we also call each matrix of this set the pattern of $\ell$. This concept was introduced 
in \cite{Gavrilyuk} where the following result was shown.

\begin{res} [\cite{Gavrilyuk}]\label{res-ti2}
Let $\mathcal{L}$ be a Cameron-Liebler line class with parameter $x$ in $\PG(3,q)$.
Let $\ell$ be a line of $\PG(3,q)$, and $T=(t_{ij})$ its pattern with respect to $\mathcal{L}$.
\renewcommand{\labelenumi}{\rm (\alph{enumi})}
\begin{enumerate}
\item For all $k,l\in\{1,\dots,q+1\}$
$$\sum_{i=1}^{q+1}t_{il}+\sum_{j=1}^{q+1}t_{kj}=
\left\{
\begin{matrix}
x+(q+1)t_{kl}&\mbox{if $\ell\notin \mathcal{L}$} &  & \\
x+(q+1)(t_{kl}+1)-2&\mbox{if $\ell\in \mathcal{L}$}.
\end{matrix}
\right.$$
\item $t_{kl}+t_{rs}=t_{ks}+t_{rl}$ for all $k,l,r,s\in\{1,\dots,q+1\}$.
\item
$$
\sum_{i,j=1}^{q+1} t_{ij}^2=
\left\{
\begin{matrix}
x(q+x)&\mbox{if $\ell\notin \mathcal{L}$} &  & \\
q^3+q^2+(x-1)^2+q(x-1)&\mbox{if $\ell\in \mathcal{L}$}.
\end{matrix}
\right.$$
\end{enumerate}
\end{res}

It is known, see \cite[Lemma~4.4]{GovaertsPenttila}, 
that a Cameron-Liebler line class $\cal L$ with parameter 
$x$ in $\PG(3,q)$ has the following property. If $P$ is a point and $\pi$ a plane of $\PG(3,q)$ with $P\in\pi$, then the number of lines of $\cal L$ through $P$ plus the number of lines of $\cal L$ in $\pi$ 
is congruent to $x$ modulo $q+1$ (cf. Result \ref{res-ti2}(a)). Thus, if $n$ is the number of lines of 
$\cal L$ in some plane, then every plane has congruent to $n$ modulo $q+1$ lines of $\cal L$, and 
the number of lines of $\cal L$ through any point is congruent to $x-n$ modulo $q+1$.
The following result determines all admissible values for $n$, 
and it provides a strong necessary existence criterion for Cameron-Liebler line classes in $\PG(3,q)$.

\begin{res} [\cite{GavrilyukMetsch}]\label{res-modular}
Suppose $\cal L$ is a Cameron-Liebler line class with parameter $x$ in $\PG(3,q)$. 
Then, for every plane and every point of $\PG(3,q)$, one has
\begin{equation}\label{eqn_main}
{x\choose 2}+n(n-x)\equiv 0 \mod q+1
\end{equation}
where $n$ is the number of lines of $\cal L$ in the plane respectively through the point.

Thus, if $\PG(3,q)$ has a Cameron-Liebler line class with parameter $x$, then
Eq. \eqref{eqn_main} has a solution for $n$ in the set $\{0,1,\dots,q\}$.
\end{res}

Result \ref{res-modular} rules out roughly at least one half of all possible values 
of $x$ from the set $\{0,1,2,\ldots,\lfloor\frac{q^2+1}{2}\rfloor\}$ 
(see \cite[Section~3]{GavrilyukMetsch}). However, it gives only a necessary but not sufficient 
criterion for the existence of Cameron-Liebler line classes: for example, 
in $\PG(3,5)$, for $x\in \{4,5,6,8,9,10,12,13\}$, Eq. (\ref{eqn_main}) 
has a solution in $n$, however, line classes with parameter $x\in \{4,5,6,8,9\}$ do not exist, 
\cite{GavrilyukMetsch}. The proof of their non-existence, which was given in \cite{GavrilyukMetsch}, 
required some ad-hoc combinatorial analysis of admissible patterns. 
In the next section, we show that this result can be obtained within a more general approach.

\section{Counting patterns}\label{section-main}

Suppose that $\PG(3,q)$ contains a Cameron-Liebler line class $\mathcal{L}$ with parameter $x$.
From the solutions of Eq. (\ref{eqn_main}) for $n$, one can determine 
all admissible patterns (see Section \ref{section-generating}), 
i.e., those square matrices of size $q+1$ which satisfy (a)-(c) of Result \ref{res-ti2}. 
Which of them can really appear as the patterns of lines with respect to $\mathcal{L}$?
How many lines of $\PG(3,q)$ may have a given pattern with respect to $\mathcal{L}$?
It turns out that non-equivalent Cameron-Liebler line classes of $\PG(3,q)$ 
with the same parameter $x$ (for example, $\mathcal{D}$, $\mathcal{R}^+$, $\mathcal{P}$ 
described in the introduction) may have different sets of patterns, and thus these questions 
cannot be answered without knowing the structure of $\mathcal{L}$.

On the other hand, it is clear that the set of patterns w.r.t. $\mathcal{L}$ 
being a subset of the set of all admissible patterns should satisfy 
certain compatibility properties. In this section, we derive some 
necessary properties, which can be expressed as a system 
of Diophantine equations (see Section \ref{section-equations}). 
Any (non-negative integer) solution of this system 
describes a possible structure of $\mathcal{L}$ (with regards to the questions above).
Conversely, if the system has no solution, this yields the non-existence 
of Cameron-Liebler line classes with parameter $x$ in $\PG(3,q)$. 
In particular, this shows in a uniform way all non-existence results for Cameron-Liebler line classes 
in $\PG(3,q)$ for $q\leq 5$ (see Section \ref{section-comp}).
Unfortunately, the existence of a non-negative integer solution 
does not necessarily imply the existence of $\mathcal{L}$, 
so further analysis may be needed (see Section \ref{sect-13-1}).

\subsection{Generating weights and patterns}\label{section-generating}

For a line $\ell$ of $\PG(3,q)$, 
each row (column) sum $s$ of the pattern of $\ell$ with respect to $\mathcal{L}$ 
is the number of lines of $\mathcal{L}\setminus \{\ell\}$ through a point (in a plane, 
respectively) which are incident to $\ell$. Then, by Result \ref{res-modular}, 
the number $n:=s+\chi_{\mathcal{L}}(\ell)$ is a solution of Eq. (\ref{eqn_main}).
We call $n$ the {\it weight} of a point (plane, respectively).
For definiteness, suppose that $s$ is the row sum. Then we define:
\begin{align*}
N&:=\{w\in \{0,1,\ldots,q^2+q+1\}\mid w\equiv n \mod q+1\},\\
M&:=\{w\in \{0,1,\ldots,q^2+q+1\}\mid w\equiv x-n \mod q+1\}.
\end{align*}

By Result \ref{res-modular} and the preceding discussion, $N$ and $M$
consist of all admissible weights of points and planes with respect to $\mathcal{L}$. 
(We note that Eq. (\ref{eqn_main}) may admit other solutions, not congruent to $n$ or $x-n$
modulo $q+1$, i.e., there may exist several sets of admissible weights.)
For a weight $w$, denote by $n_w$ (by $m_w$) the number of points (planes, respectively) of weight $w$.

\begin{lemma}\label{lemma-planes-counting}
The following holds:
\begin{align}
\label{eq-sum-n}
\sum_{u\in N} n_u = \sum_{w\in M} m_w &= q^3+q^2+q+1,\\
\label{eq-sum-wn}
\sum_{u\in N} un_u = \sum_{w\in M} wm_w &= x(q^2+q+1)(q+1),\\
\label{eq-sumww1n}
\sum_{u\in N} u(u-1)n_u = \sum_{w\in M} w(w-1)m_w &= x(q^2 + q + 1)((q + 1)x + q^2 - 1).
\end{align}
\end{lemma}
\proof 
Eq. (\ref{eq-sum-n}) is straightforward (its right-hand side equals the number of points/planes).
Recall that, by Result \ref{res-basic}(a), $|\mathcal{L}|=x(q^2+q+1)$ holds. 
Every line lies in $q+1$ planes and consists of $q+1$ points, which gives Eq. (\ref{eq-sum-wn}).
For Eq. (\ref{eq-sumww1n}), we count the number of pairs of intersecting lines of $\mathcal{L}$: 
by Result \ref{res-basic}(b), each line of $\mathcal{L}$ intersects $(q + 1)x + q^2 - 1$ 
other lines of $\mathcal{L}$, while every two distinct lines in a plane (on a point) 
have a point in common.
\wbull

Observe that, by Result \ref{res-ti2}(b), each pattern is determined by 
any of its row-column pairs. 
Let $T_1,\ldots,T_{k_1}$ be all admissible patterns for lines of $\mathcal{L}$, 
and $T_{k_1+1},\ldots,T_{k_1+k_0}$ all admissible patterns for lines of $\overline{\mathcal{L}}$.
Let $I_1:=\{1,2\ldots,k_1\}$ and 
$I_0:=\{k_1+1,\ldots,k_1+k_0\}$ denote the corresponding sets of indices of patterns, 
and set $I:=I_0\cup I_1$.


%

\subsection{A system of equations}\label{section-equations}

For a pattern $T_i$, $i\in I$, denote by $c_{i,w}$ (by $r_{i,w}$) the numbers of its columns 
(its rows, resp.) summing to $w-\chi_{I_1}(i)$.
We now introduce our main unknowns: let $z_i$, $i\in I$, denote the number 
of lines (from ${\cal L}$ if $i\in I_1$, or from $\overline{\cal L}$ if $i\in I_0$) 
having pattern $T_i$.

\begin{lemma}
The following holds.

\begin{equation}\label{eq-allsum}
\begin{matrix}
{\displaystyle\sum_{i\in I_1}z_i=|{\cal L}|=x(q^2+q+1),} & 
{\displaystyle\sum_{i\in I_0}z_i=|\overline{{\cal L}}|=(q^2+1-x)(q^2+q+1).} \\
\end{matrix}
\end{equation}
For all $w,w'\in M$ with $w\ne w'$ and all $u,u'\in N$ with $u\ne u'$:
\begin{equation}\label{eq-colsum}
\begin{matrix}
{\displaystyle\sum_{i\in I}c_{i,w}z_i=(q^2+q+1)m_w,} & 
{\displaystyle\sum_{i\in I}r_{i,u}z_i=(q^2+q+1)n_{u},}\\
\end{matrix}
\end{equation}
\begin{equation}\label{eq-wcolsum}
\begin{matrix}
{\displaystyle\sum_{i\in I_1}c_{i,w}z_i=wm_w,} & 
{\displaystyle\sum_{i\in I_1}r_{i,u}z_i=wn_{u},} \\
\end{matrix}
\end{equation}
\begin{equation}\label{eq-pairsum-1}
\begin{matrix}
{\displaystyle\sum_{i\in I}c_{i,w}(c_{i,w}-1)z_i=m_w(m_w-1),} & 
{\displaystyle\sum_{i\in I}r_{i,u}(r_{i,u}-1)z_i=n_u(n_u-1),} \\
\end{matrix}
\end{equation}
\begin{equation}\label{eq-pairsum-2}
\begin{matrix}
{\displaystyle\sum_{i\in I}c_{i,w}c_{i,w'}z_i=m_wm_{w'},} & 
{\displaystyle\sum_{i\in I}r_{i,u}r_{i,u'}z_i=n_un_{u'},} \\
\end{matrix}
\end{equation}
\begin{equation}\label{eq-pencilscounting}
{\displaystyle\Big(q+1-\frac{w+u-x}{q+1}\Big)\sum_{i\in I_1} c_{i,w}r_{i,u}z_i = 
\frac{w+u-x}{q+1}\sum_{j\in I_0} c_{j,w}r_{j,u}z_j.}
\end{equation}
\end{lemma}
\proof Eq. (\ref{eq-allsum}) is straightforward. 
Each of Eqs. (\ref{eq-colsum})-(\ref{eq-pairsum-2}) is obtained by 
double counting, and we prove only its left part, as the right one follows by a point-plane duality.
For Eq. (\ref{eq-colsum}) (Eq. (\ref{eq-wcolsum})) 
we count the number of pairs $(\ell,\pi)$ where $\ell$ is a line (a line of $\mathcal{L}$, respectively)
in a plane $\pi$ of weight $w$.
For Eq. (\ref{eq-pairsum-1}) (Eq. (\ref{eq-pairsum-2})) we count the number of pairs of 
different planes of the same weight $w$ (of different weights $w,w'$, respectively).

Finally, in order to obtain Eq. (\ref{eq-pencilscounting}) we count the number of 
pairs of concurrent lines $(\ell_1,\ell_2)$ with $\ell_1\in \mathcal{L}$, $\ell_2\in\overline{\mathcal{L}}$ 
and $\ell_1,\ell_2\in \St(P)\cap \Li(\pi)$ where a plane $\pi$ has weight $w$ and a point $P$ 
has weight $u$.
By Result \ref{res-ti2}(a), such a pencil $\St(P)\cap \Li(\pi)$ contains precisely 
$(w+u-x)/(q+1)$ lines from $\mathcal{L}$. The lemma is proved.
\wbull

\subsection{Computational aspects and results}\label{section-comp}

In this section, we discuss solving the system of (Diophantine) equations 
(\ref{eq-sum-n})-(\ref{eq-pencilscounting}) with respect to the unknowns 
$m_w$ for all $w\in M$, $n_u$ for all $u\in N$, and $z_i$ for all $i\in I$.
Note that this system is not linear with respect to $m_w$ and $n_u$.
The system of equations (\ref{eq-allsum})-(\ref{eq-pencilscounting}) 
is linear with respect to the $|I|$ unknowns $z_i$, and it consists of at most 
$2+4(|M|+|N|)+{|M|\choose 2}+{|N|\choose 2}+|M||N|$ equations, however, 
it appears to be hard to determine its rank (even calculating 
the precise values of $|M|$, $|N|$ can be very tedious, see \cite[Section~3]{GavrilyukMetsch}), 
as its coefficients depend on the structure of all patterns.

We consider the following approach to solving Eqs. (\ref{eq-sum-n})-(\ref{eq-pencilscounting}).
First, we regard the system of Eqs. (\ref{eq-allsum})-(\ref{eq-pencilscounting}) 
as a system of linear equations with respect to the unknowns $z_i$ (so that its left-hand side 
consists of linear combinations of $z_i$ only), and apply the Gaussian elimination 
procedure to it, in which $n_u$, $m_w$ 
are treated as indeterminates. 
It may happen that the row-reduced echelon form of the system 
contains {\it zero} equations in its left-hand side with a non-zero right-hand side, which 
is represented by polynomials in $n_u$, $m_w$ of degree at most $2$.
Therefore, in order the system of Eqs. (\ref{eq-allsum})-(\ref{eq-pencilscounting}) 
to be consistent the polynomials in the right-hand side of the zero equations 
must be equated to $0$. This may provide additional constraints on a solution of 
Eqs. (\ref{eq-sum-n})-(\ref{eq-sumww1n}).
Using these additional constraints, with the aid of computer we find all solutions 
of Eqs. (\ref{eq-sum-n})-(\ref{eq-sumww1n}) (perhaps, one can use a more sophisticated 
approach, e.g., based on a Gr\"obner basis, however, for our examples considered below 
a brute force search was sufficient). 
Let us call such a solution 
of Eqs. (\ref{eq-sum-n})-(\ref{eq-sumww1n}) a {\it feasible weight distribution} of planes and points.
Further, we substitute every feasible weight distribution back into 
the row-reduced echelon form of the system of 
Eqs. (\ref{eq-allsum})-(\ref{eq-pencilscounting}) and so obtain 
a system of linear Diophantine equations with respect to the unknowns $z_i$, 
and finally solve it with the aid of the MIP solver Gurobi \cite{Gurobi} 
(in general, its rank is less than the number of unknowns).

We now summarize the results for non-trivial Cameron-Liebler line classes in $\PG(3,q)$ with $q\leq 5$.
In $\PG(3,2)$, all line classes are trivial.
In $\PG(3,3)$, for $x\in \{3,4\}$, Eq. (\ref{eqn_main}) has no solution, and for $x=5$ 
there exists a line class 
(the first counter example to the Cameron-Liebler conjecture found by Drudge \cite{Drudge}), 
which is the only non-trivial Cameron-Liebler line class in $\PG(3,3)$. 
Let us consider this example in detail in order to illustrate the equations from 
Sections \ref{section-generating} and \ref{section-equations}.

For $q=3$ and $x=5$, Eq. (\ref{eqn_main}) admits only two sets of (pairwise congruent modulo $q+1$) 
solutions:
\[
N=\{2,6,10\},~~M=\{3,7,11\}, 
\]
which, by Result \ref{res-ti2}, give the following admissible patterns for lines of $\mathcal{L}$:
\[\small
T_1=
\begin{pmatrix}
0& 0& 0& 1\cr
2& 2& 2& 3\cr
2& 2& 2& 3\cr
2& 2& 2& 3
\end{pmatrix}, 
T_2=
\begin{pmatrix}
0& 1& 2& 2\cr
0& 1& 2& 2\cr
1& 2& 3& 3\cr
1& 2& 3& 3
\end{pmatrix},
\]
and for lines of $\overline{\mathcal{L}}$:
\[\small
T_3=
\begin{pmatrix}
0& 0& 0& 2\cr
1& 1& 1& 3\cr
1& 1& 1& 3\cr
1& 1& 1& 3
\end{pmatrix}, 
T_4=
\begin{pmatrix}
0& 0& 1& 1\cr
0& 0& 1& 1\cr
1& 1& 2& 2\cr
2& 2& 3& 3
\end{pmatrix}.
\]

Further, the system of Eqs. (\ref{eq-sum-n})-(\ref{eq-sumww1n}):
\begin{align*}
n_2+n_6+n_{10}&=40 & m_3+m_7+m_{11}&=40,\\
2n_2+6n_6+10n_{10}&=260 & 3m_3+7m_7+11m_{11}&=260,\\
2n_2+30n_6+90n_{10}&=1820 & 6m_3+42m_7+110m_{11}&=1820 
\end{align*}
has a unique solution:
\begin{align*}
n_2&=10,~n_6=15,~n_{10}=15, & m_3&=15,~m_7=15,~m_{11}=10.
\end{align*}

The system of Eqs. (\ref{eq-allsum})-(\ref{eq-pencilscounting}) consists of $33$ equations, 
its row-reduced echelon form (with respect to the unknowns $z_1$, $z_2$, $z_3$, $z_4$)
contains 4 non-zero equations:
\begin{align*}
z_1  ~~   ~~   ~~  &=20,\\
~~   z_2  ~~   ~~  &=45,\\
~~   ~~   z_3  ~~  &=20,\\
~~   ~~   ~~   z_4 &=45,
\end{align*}
and 18 zero equations with a non-zero right-hand side which have form:
\[
0=y+C \mbox{~or~} 0=yy'+C,
\]
where $y,y'\in \{n_2,n_6,n_{10}\}$ or $y,y'\in \{m_3,m_7,m_{11}\}$, and 
$C$ is a number.
Clearly, for this example, these zero equations are redundant, as 
the only feasible weight distribution is uniquely determined from Eqs. (\ref{eq-sum-n})-(\ref{eq-sumww1n}).

We obtain that $\mathcal{L}$ 
consists of $20$ lines with pattern $T_1$ and $45$ lines with pattern $T_2$, 
and its complementary set $\overline{\mathcal{L}}$ consists of 
$20$ lines with pattern $T_3$ and $45$ lines with pattern $T_4$. 
One can see that the set $\mathcal{P}_2$ of points of weight $2$ is a $10$-{\it cap} 
(and thus it is an elliptic quadric, see \cite{Barlotti,Panella}), 
i.e., every line of $\PG(3,3)$ intersects $\mathcal{P}_2$ in at most two points, 
while the $45={10\choose 2}$ lines with pattern $T_4$ are the 2-secants to $\mathcal{P}_2$.
It is not difficult to further recover the structure of $\mathcal{L}$ as it is described 
in \cite{Drudge}, and to prove its uniqueness.

For $\PG(3,4)$, Eq. (\ref{eqn_main}) has no solution when $x\in \{3,4,8\}$. 
Further, we obtain:
\begin{itemize}
\item if $x\in \{5,6\}$, then Eqs. (\ref{eq-sum-n})-(\ref{eq-sumww1n}) together with 
zero equations of the row-reduced echelon form of Eqs. (\ref{eq-allsum})-(\ref{eq-pencilscounting})
have no solution in non-negative integers. Thus, no such Cameron-Liebler line class exists.

\item if $x=7$, then $N=M=\{1,6,11,16,21\}$. 
We have 7 admissible patterns for lines of $\mathcal{L}$, 
and 4 for lines of $\overline{\mathcal{L}}$. The system of 
Eqs. (\ref{eq-allsum})-(\ref{eq-pencilscounting}) consists of 62 equations, 
while its row-reduced echelon form consists of 49 equations 
and has rank 11 with respect to $z_i$ (so applying the Gaussian elimination procedure 
gives $49-11=38$ zero equations with respect to $z_i$).
The system of Eqs. (\ref{eq-sum-n})-(\ref{eq-sumww1n}) together with these 38 equations 
gives a unique feasible weight distribution, and the row-reduced echelon form 
of Eqs. (\ref{eq-allsum})-(\ref{eq-pencilscounting}) admits a unique solution 
for $z_i$. 
The corresponding Cameron-Liebler line class does exist, it was the first example 
in even characteristic found by Govaerts and Penttila in \cite{GovaertsPenttila}, 
its uniqueness was shown in \cite{Gavrilyuk}.
\end{itemize}

For $\PG(3,5)$, Eq. (\ref{eqn_main}) has no solution when $x\in \{3,7,11\}$.
Further, we obtain:
\begin{itemize}
\item if $x=4$, then Eq. (\ref{eqn_main}) has a solution in $n$, but 
there are no admissible patterns satisfying Result \ref{res-ti2}(c).
If $x\in \{5,6,8,9\}$, then Eqs. (\ref{eq-sum-n})-(\ref{eq-sumww1n}) together with 
zero equations of the row-reduced echelon form of Eqs. (\ref{eq-allsum})-(\ref{eq-pencilscounting})
have no solution in non-negative integers. Thus, no such Cameron-Liebler line class exists.

\item if $x=10$, then $N=\{1, 7, 13, 19, 25, 31\}$, 
$M=\{3, 9, 15, 21, 27\}$.
We have 10 admissible patterns for lines of $\mathcal{L}$, 
and 6 for lines of $\overline{\mathcal{L}}$. The system of 
Eqs. (\ref{eq-allsum})-(\ref{eq-pencilscounting}) consists of 76 equations, 
while its row-reduced echelon form consists of 62 equations 
and has rank 16 with respect to $z_i$ (so applying the Gaussian elimination procedure 
gives $62-16=48$ zero equations with respect to $z_i$).
The system of Eqs. (\ref{eq-sum-n})-(\ref{eq-sumww1n}) together with these 48 equations 
gives a unique feasible weight distribution, and the row-reduced echelon form 
of Eqs. (\ref{eq-allsum})-(\ref{eq-pencilscounting}) admits a unique solution for $z_i$. 
The corresponding Cameron-Liebler line class $\mathcal{G}$ does exist, 
it was constructed and shown to be unique in \cite{GavrilyukMetsch}.

\item if $x=12$, then $N=M=\{0, 6, 12, 18, 24, 30\}$. 
We have 8 admissible patterns for lines of $\mathcal{L}$, 
and 8 for lines of $\overline{\mathcal{L}}$. The system of 
Eqs. (\ref{eq-allsum})-(\ref{eq-pencilscounting}) consists of 84 equations, 
while its row-reduced echelon form consists of 69 equations 
and has rank 16 with respect to $z_i$ (so applying the Gaussian elimination procedure 
gives $69-16=53$ zero equations with respect to $z_i$).
The system of Eqs. (\ref{eq-sum-n})-(\ref{eq-sumww1n}) together with these 53 equations 
gives a unique feasible weight distribution, and the row-reduced echelon form 
of Eqs. (\ref{eq-allsum})-(\ref{eq-pencilscounting}) admits a unique solution for $z_i$. 
As we mentioned in the introduction, the Cameron-Liebler line class $\mathcal{R}$, which 
was constructed in \cite{Rodgers} and generalised to an infinite family in \cite{Metsch-and-Co} and 
\cite{FMX}, has parameter $x=12$, however, its uniqueness was not known. 
In Section \ref{sect-12}, we analyse the solution of Eqs. (\ref{eq-sum-n})-(\ref{eq-pencilscounting}) 
and prove that $\mathcal{R}$ is the only Cameron-Liebler line class with parameter $12$ in $\PG(3,5)$.

\item if $x=13$, then we have the following two cases: 
$N=\{0, 6, 12, 18, 24, 30\}$, $M=\{1, 7, 13, 19, 25, 31\}$ 
or $N=\{3, 9, 15, 21, 27\}$, $M=\{4, 10, 16, 22, 28\}$.
In each of them, there are 10 admissible patterns for lines of $\mathcal{L}$, 
and 10 for lines of $\overline{\mathcal{L}}$.

In the first case, the system of 
Eqs. (\ref{eq-allsum})-(\ref{eq-pencilscounting}) consists of 84 equations, 
while its row-reduced echelon form consists of 69 equations 
and has rank 20 with respect to $z_i$ (so applying the Gaussian elimination procedure 
gives $69-20=49$ zero equations with respect to $z_i$).
The system of Eqs. (\ref{eq-sum-n})-(\ref{eq-sumww1n}) together with these 49 equations 
gives 3 feasible weight distributions. Substituting each of them into the 
row-reduced echelon form of Eqs. (\ref{eq-allsum})-(\ref{eq-pencilscounting}), 
we obtain a unique solution for $z_i$, thus, in total the three solutions of 
Eqs. (\ref{eq-sum-n})-(\ref{eq-pencilscounting}). 
These solutions will be analysed in Section \ref{sect-13-1}.
We note that one of them corresponds to the Cameron-Liebler line class $\mathcal{R}^+$.

In the second case, the system of 
Eqs. (\ref{eq-allsum})-(\ref{eq-pencilscounting}) consists of 71 equations, 
while its row-reduced echelon form consists of 57 equations 
and has rank 20 with respect to $z_i$ (so applying the Gaussian elimination procedure 
gives $57-20=37$ zero equations with respect to $z_i$).
The system of Eqs. (\ref{eq-sum-n})-(\ref{eq-sumww1n}) together with these 37 equations 
gives 2 feasible weight distributions, and, for both of them, 
the row-reduced echelon form of Eqs. (\ref{eq-allsum})-(\ref{eq-pencilscounting}) 
admits a unique solution for $z_i$. 
These solutions correspond to the Cameron-Liebler line classes $\mathcal{P}$ and $\mathcal{D}$, 
see Section \ref{sect-13-2}.
\end{itemize}

\section{Cameron-Liebler line classes in $\PG(3,5)$}\label{section-5}

In this section, we prove our main result, Theorem \ref{theo-main}.
Let $\mathcal{L}$ be a non-trivial Cameron-Liebler line class with parameter $x$ in $\PG(3,5)$.
It follows from Section \ref{section-comp} that $x\in \{10,12,13\}$, 
and if $x=10$, then $\mathcal{L}$ or its dual line class is projectively equivalent 
to $\mathcal{G}$.
If $x=13$, then Eqs. (\ref{eq-sum-n})-(\ref{eq-pencilscounting}) admit 5 solutions, 
which can be divided into the two groups according to the sets of possible weights: 
$N=\{0, 6, 12, 18, 24, 30\}$, $M=\{1, 7, 13, 19, 25, 31\}$ (three solutions)
or $N=\{3, 9, 15, 21, 27\}$, $M=\{4, 10, 16, 22, 28\}$ (two solutions).
In Section \ref{sect-13-1}, we rule out two solutions from the first group, 
and observe that the only solution, which is left, corresponds to a Cameron-Liebler 
line class that can be seen as a disjoint union of lines of a Cameron-Liebler line 
class $\mathcal{L}'$ with parameter $x=12$ 
and the lines of some plane disjoint from $\mathcal{L}'$.
Thus, the question in this case reduces to $x=12$, 
and by analysing the solution of Eqs. (\ref{eq-sum-n})-(\ref{eq-pencilscounting}) for $x=12$, 
we show in Section \ref{sect-12} that $\mathcal{L}'$ is projectively equivalent to $\mathcal{R}$, 
and hence $\mathcal{L}$ is projectively equivalent to $\mathcal{R}^+$.
In Section \ref{sect-13-2}, we show that the two solutions of the second group 
correspond to the Cameron-Liebler line classes $\mathcal{D}$ and $\mathcal{P}$ only.

\subsection{The first group of solutions in the case $x=13$}\label{sect-13-1}

The patterns of lines, the feasible weight distributions and the solutions of 
Eqs. (\ref{eq-allsum})-(\ref{eq-pencilscounting}) corresponding to 
the case $N=\{0, 6, 12, 18, 24, 30\}$, $M=\{1, 7, 13, 19, 25, 31\}$ 
are given in Tables \ref{patterns-13-1-1}, \ref{patterns-13-1-0}, 
\ref{table-13-1}, \ref{table-13-2}, respectively.
(Tables \ref{patterns-13-1-1}-\ref{table-13-2-2} are given in the appendix).
We will refer to them as the solutions $\#1$, $\#2$, and $\#3$ of the first group.

\begin{lemma}\label{lemma-group1}
There are no Cameron-Liebler line classes corresponding to the solutions $\#1$ or $\#2$.
\end{lemma}
\proof On the contrary, suppose that a Cameron-Liebler line class $\mathcal{L}$ with 
parameter $x=13$ in $\PG(3,5)$ corresponds to the solution $\#1$ or $\#2$. 
It follows from Table \ref{table-13-1} that there exists 
a point $P$ with $|\St(P)\cap \mathcal{L}|=12$. Consider the quotient geometry $\PG(3,5)/P$ 
on the point $P$, i.e., a projective plane $\pi\cong \PG(2,5)$ whose points are the lines on $P$ 
and whose lines are the planes on $P$. The set $\St(P)\cap \mathcal{L}$ corresponds to a set $B$ 
of 12 points of $\pi$ such that, for any point $p\in \pi$, 
the row $i_p:=(|\ell\cap (B\setminus \{p\})|: \ell\in \St(p))$ is permutation-equivalent 
to one of the following:
\[\small
t_2=
\begin{pmatrix}
1&  1&  2&  2&  2&  3
\end{pmatrix},
t_3=
\begin{pmatrix}
0&  1&  1&  3&  3&  3
\end{pmatrix},
t_4=
\begin{pmatrix}
0&  1&  2&  2&  3&  3
\end{pmatrix},
\]
\[\small
t_5=
\begin{pmatrix}
0&  2&  2&  2&  2&  3
\end{pmatrix},
t_6=
\begin{pmatrix}
1&  1&  1&  2&  3&  3
\end{pmatrix},
\]
if $p\in B$, or 
\[\small
t_{10}=
\begin{pmatrix}
1& 1& 1& 2& 3& 4
\end{pmatrix},
t_{11}=
\begin{pmatrix}
1& 1& 1& 3& 3& 3
\end{pmatrix},
t_{13}=
\begin{pmatrix}
1& 1& 2& 2& 2& 4
\end{pmatrix},\]
\[\small
t_{14}=
\begin{pmatrix}
1& 1& 2& 2& 3& 3
\end{pmatrix},
t_{16}=
\begin{pmatrix}
2& 2& 2& 2& 2& 2
\end{pmatrix},
\]
if $p\notin B$ (so that $t_i$ is the row of $T_i$ that corresponds to a point of weight $12$).

Moreover, it follows from Table \ref{table-13-2} that $z_5=1$, and so 
we may assume that there exists a line $\ell^{\star}$ on $P$ with pattern $T_5$, 
and hence there exists a unique point $p^{\star}\in B$ with $i_{p^{\star}}=t_5$.
An exhaustive computer search shows that no such set $B$ exists in $\PG(2,5)$. 
\wbull

Let $\mathcal{L}$ be a Cameron-Liebler line class with parameter $x=13$ in $\PG(3,5)$ 
that satisfies the solution $\#3$. It follows from Table \ref{table-13-1} that there exists 
a unique plane $\pi_0$ with $|\Li(\pi_0)\cap \mathcal{L}|=31$.
Now, $\mathcal{L}':=\mathcal{L}\setminus \Li(\pi_0)$ is a Cameron-Liebler line class 
with parameter $x=12$ in $\PG(3,5)$, and we consider this line class in the next section.

\subsection{A Cameron-Liebler line class with parameter $x=12$}\label{sect-12}

In the case $x=12$, the sets of possible weights are $N=M=\{0, 6, 12, 18, 24, 30\}$, 
and there exists the only feasible weight distribution:
\[
n_0=1, n_6=31, n_{12}=62, n_{18}=31, n_{24}=31, n_{30}=0,
\]
and $m_w=n_w$ for all $w\in M$. Further, the patterns of lines are given in Tables 
\ref{patterns-12-1}, \ref{patterns-12-0}, and the only solution of 
Eqs. (\ref{eq-allsum})-(\ref{eq-pencilscounting}) is given in Table \ref{table-12}, 
respectively.

\begin{lemma}\label{lemma-12}
$\mathcal{L}'$ is projectively equivalent to $\mathcal{R}$.
\end{lemma}
\proof
In order to prove the uniqueness of $\mathcal{R}$, we make use of the following 
observation, which is based on Result \ref{res-basic}(c). For a fixed line $\ell^{\star}$ 
of $\PG(3,q)$, suppose we are given the set $\mathcal{L}'_1$ of lines of $\mathcal{L}'$ 
that intersect $\ell^{\star}$. Then one can uniquely determine the set $\mathcal{L}'\setminus \mathcal{L}'_1$
by evaluating, for every line $\ell$ skew to $\ell^{\star}$, the number of lines of $\mathcal{L}'_1$ 
that intersect both $\ell^{\star}$ and $\ell$.
Thus, all Cameron-Liebler line classes with given parameter $x$ in $\PG(3,q)$
can be found by applying this procedure to all possible sets $\mathcal{L}_1'$ 
that are compatible with the pattern of line $\ell^{\star}$.

Let $e_1,\ldots,e_4$ denote the projective points of $\PG(3,5)$ 
corresponding to the standard basis of $\mathbb{F}_5^4$.
For our purpose, we fix a line $\ell^{\star}:=\langle e_1,e_2\rangle$ and suppose that 
it has pattern $T_6$ represented in Table \ref{table-t6}, 
where, without loss of generality, $\pi_1:=\langle e_1,e_2,e_3\rangle$, 
$\pi_2:=\langle e_1,e_2,e_4\rangle$, $\pi_3:=\langle e_1,e_2,e_3+e_4\rangle$.
Our approach is as follows. First, we determine all possible planar sections of $\mathcal{L}'$.
Note that $\mathcal{L}'$ is formally self-dual in the sense that its dual line class 
has the same set of patterns and also satisfies the solution in Table \ref{table-12}.
Further, one can consider all possible realizations of the planar sections in all planes on line $\ell^{\star}$, which thereby determine all candidates for $\mathcal{L}'_1$.
In order to reduce the amount of computations, we only consider 
planar sections in the planes $\pi_1$, $\pi_2$, and $\pi_3$ 
(with two cases for the weight of $\pi_3$, which are marked in Table \ref{table-t6}).
For given planar sections in $\pi_1$, $\pi_2$, and $\pi_3$, we then consider 
the quotient geometries on the points of $\ell^{\star}$. 
Note that the planar sections in $\pi_1$, $\pi_2$, and $\pi_3$ 
already determine all the lines of $\mathcal{L}'$ on $e_1$.
For $e_2$, these planar sections determine $3+4+5=12$ lines out of 18 lines 
of $\mathcal{L}'$ on $e_2$, and we consider all planar sections of weight 18 
corresponding to the quotient geometry on $e_2$ that match these 12 lines.
This gives all candidates for the lines of $\mathcal{L}'$ on $e_2$.
The same procedure applies to all points of $\ell^{\star}$, 
which finally gives all candidates for $\mathcal{L}'_1$.

Let us describe some details.
An analysis similar to that of the proof of Lemma \ref{lemma-group1}
shows that the planar sections of $\mathcal{L}'$ are projectively 
equivalent to the following configurations of lines in $\PG(2,5)$:
\begin{itemize}
\item weight $6$: $\mathcal{S}_6:=\{\langle 1:0:0\rangle$,
$\langle 0:0:1\rangle$,
$\langle 0:1:2\rangle$,
$\langle 1:3:3\rangle$,
$\langle 1:3:0\rangle$,
$\langle 1:4:0\rangle\}$.

In fact, an analysis similar to that of the proof of Lemma \ref{lemma-group1} 
gives two non-equivalent configurations of six lines in $\PG(2,5)$. 
One of them, namely $\mathcal{S}_6$, consists of three lines that should have pattern $T_3$ 
and three lines that should have pattern $T_4$, 
while the other one consists of six lines that may have pattern $T_4$ only.
As there exist $m_6=31$ planes of weight $6$ with respect to $\mathcal{L}'$, and 
two such planes cannot share a line of $\mathcal{L}'$, it follows that the union of these 31 planes 
contains $31\cdot 6=186$ lines of $\mathcal{L}'$, while $z_3=z_4=93$. 
Therefore, the second configuration cannot appear as a planar section with respect 
to $\mathcal{L}'$. 

\item weight $12$:
$\mathcal{S}^{(1)}_{12}:=\{\langle 1 : 0 : 0  \rangle$, 
$\langle 0 : 1 : 2  \rangle$, 
$\langle 1 : 4 : 3  \rangle$, 
$\langle 0 : 1 : 1  \rangle$, 
$\langle 1 : 4 : 4  \rangle$, 
$\langle 0 : 0 : 1  \rangle$, 
$\langle 1 : 3 : 3  \rangle$, 
$\langle 0 : 1 : 0  \rangle$, 
$\langle 1 : 4 : 0  \rangle$, 
$\langle 1 : 0 : 3  \rangle$, 
$\langle 1 : 2 : 0  \rangle$, 
$\langle 1 : 2 : 1  \rangle\}$, 
or 
$\mathcal{S}^{(2)}_{12}:=\{\langle 1 : 0 : 0  \rangle$, 
$\langle 1 : 1 : 1  \rangle$, 
$\langle 1 : 2 : 0  \rangle$, 
$\langle 1 : 2 : 2  \rangle$, 
$\langle 1 : 3 : 0  \rangle$, 
$\langle 1 : 3 : 4  \rangle$, 
$\langle 1 : 0 : 3  \rangle$, 
$\langle 1 : 0 : 1  \rangle$, 
$\langle 1 : 3 : 2  \rangle$, 
$\langle 1 : 1 : 4  \rangle$, 
$\langle 1 : 4 : 1  \rangle$, 
$\langle 1 : 3 : 1  \rangle\}$.
We note that if a planar section of weight 12 contains 
a line with pattern $T_4$, or $T_5$, or $T_6$, then it must be projectively 
equivalent to $\mathcal{S}^{(2)}_{12}$, as $\mathcal{S}^{(1)}_{12}$ 
does not contain a line with any of these patterns.

\item weight $18$: $\mathcal{S}_{18}:=\{\langle 1 : 0 : 0  \rangle$,  
$\langle 1 : 3 : 1  \rangle$,  
$\langle 1 : 4 : 0  \rangle$, 
$\langle 1 : 3 : 3  \rangle$, 
$\langle 1 : 0 : 2  \rangle$, 
$\langle 0 : 0 : 1  \rangle$, 
$\langle 0 : 1 : 1  \rangle$, 
$\langle 0 : 1 : 3  \rangle$, 
$\langle 1 : 4 : 1  \rangle$, 
$\langle 1 : 0 : 1  \rangle$, 
$\langle 1 : 2 : 1  \rangle$, 
$\langle 1 : 3 : 4  \rangle$, 
$\langle 1 : 0 : 4  \rangle$, 
$\langle 1 : 3 : 0  \rangle$, 
$\langle 0 : 1 : 2  \rangle$, 
$\langle 0 : 1 : 0  \rangle$, 
$\langle 1 : 1 : 3  \rangle$, 
$\langle 1 : 1 : 2  \rangle\}$.

\item weight $24$: $\mathcal{S}_{24}$ is the complement to 
$\{\langle 0 : 1 : 4 \rangle$, 
$\langle 1 : 0 : 4 \rangle$, 
$\langle 1 : 3 : 0 \rangle$, 
$\langle 1 : 2 : 0 \rangle$, 
$\langle 1 : 0 : 2 \rangle$, 
$\langle 1 : 0 : 1 \rangle$, 
$\langle 1 : 3 : 1 \rangle\}$.

Similarly to the case of weight $6$, there exist two admissible configurations of 24 lines in $\PG(2,5)$.  
One of them, namely $\mathcal{S}_{24}$, contains only one line that should have pattern $T_{10}$, 
while the other one contains two lines that are supposed to have pattern $T_{10}$.
The planes of weight $24$ cannot share a line with pattern $T_{10}$, 
and there exist $m_{24}=31$ planes of weight $24$ with respect to $\mathcal{L}'$ and 
$z_{10}=31$ lines with pattern $T_{10}$. Thus, every plane of weight $24$ has only 
one line with pattern $T_{10}$, and the second configuration cannot appear as a planar 
section with respect to $\mathcal{L}'$. 
\end{itemize}

Further, the set-stabiliser of $\mathcal{S}^{(2)}_{12}$ has the only orbit 
on lines that should have pattern $T_6$, so, in the plane $\pi_1$, 
we can consider an arbitrary configuration of 12 lines that is projectively equivalent 
to $\mathcal{S}^{(2)}_{12}$ and contains the line $\ell^{\star}$ 
(as it is represented in Table \ref{table-t6}).
The set-stabiliser of $\mathcal{S}_{18}$ has the only orbit 
on lines that should have pattern $T_6$. 
In the plane $\pi_2$, we consider all configurations of $18$ lines that are projectively equivalent 
to $\mathcal{S}_{18}$ and contain the line $\ell^{\star}$ 
as it is represented in Table \ref{table-t6}, and that are non-equivalent under 
the action of the point-wise stabiliser of $\ell^{\star}$ in $\PG(2,5)$.
This gives 4 candidates for a planar section of $\mathcal{L}'$ in $\pi_2$.
In the plane $\pi_3$, we consider all configurations of $6$ or $24$ lines 
that are projectively equivalent to $\mathcal{S}_{6}$ or $\mathcal{S}_{24}$, respectively, 
and contain the line $\ell^{\star}$ as it is represented in Table \ref{table-t6}.
This gives 1200 candidates in the case of weight 6 and 
400 candidates in the case of weight 24 for a planar section of $\mathcal{L}'$ in $\pi_3$.

For all these $1\times 4\times 1200 + 1\times 4\times 400$ candidates 
for the planar sections in $\pi_1,\pi_2,\pi_3$, 
the procedure as described above was completed successfully in one case only
and returned a line class $\mathcal{L}'$, which is projectively equivalent to $\mathcal{R}$.
The procedure required a few minutes of CPU time on a single core of Intel Core i5-3360M 2.8 GHz processor.
\wbull

\begin{table}[]
\centering
\caption{Pattern $T_6$ for line $\ell^{\star}$}
\label{table-t6}
\begin{tabular}{ccccccc}
                           & $\pi_1$                & $\pi_2$                & \multicolumn{2}{c}{$\pi_3$}                     &                        &                        \\
\multicolumn{1}{l}{}       & \multicolumn{1}{l}{}   & \multicolumn{1}{l}{}   & \multicolumn{2}{l}{$\swarrow\searrow$}          & \multicolumn{1}{l}{}   & \multicolumn{1}{l}{}   \\ \cline{2-7} 
\multicolumn{1}{c|}{$e_1$} & \multicolumn{1}{c|}{1} & \multicolumn{1}{c|}{2} & \multicolumn{1}{c|}{3} & \multicolumn{1}{c|}{0} & \multicolumn{1}{c|}{0} & \multicolumn{1}{c|}{0} \\ \cline{2-7} 
\multicolumn{1}{c|}{$e_2$} & \multicolumn{1}{c|}{3} & \multicolumn{1}{c|}{4} & \multicolumn{1}{c|}{5} & \multicolumn{1}{c|}{2} & \multicolumn{1}{c|}{2} & \multicolumn{1}{c|}{2} \\ \cline{2-7} 
\multicolumn{1}{c|}{}      & \multicolumn{1}{c|}{2} & \multicolumn{1}{c|}{3} & \multicolumn{1}{c|}{4} & \multicolumn{1}{c|}{1} & \multicolumn{1}{c|}{1} & \multicolumn{1}{c|}{1} \\ \cline{2-7} 
\multicolumn{1}{c|}{}      & \multicolumn{1}{c|}{2} & \multicolumn{1}{c|}{3} & \multicolumn{1}{c|}{4} & \multicolumn{1}{c|}{1} & \multicolumn{1}{c|}{1} & \multicolumn{1}{c|}{1} \\ \cline{2-7} 
\multicolumn{1}{c|}{}      & \multicolumn{1}{c|}{2} & \multicolumn{1}{c|}{3} & \multicolumn{1}{c|}{4} & \multicolumn{1}{c|}{1} & \multicolumn{1}{c|}{1} & \multicolumn{1}{c|}{1} \\ \cline{2-7} 
\multicolumn{1}{c|}{}      & \multicolumn{1}{c|}{2} & \multicolumn{1}{c|}{3} & \multicolumn{1}{c|}{4} & \multicolumn{1}{c|}{1} & \multicolumn{1}{c|}{1} & \multicolumn{1}{c|}{1} \\ \cline{2-7} 
\end{tabular}
\end{table}

\subsection{The second group of solutions in the case $x=13$}\label{sect-13-2}

The patterns of lines, the feasible weight distributions and the solutions of 
Eqs. (\ref{eq-allsum})-(\ref{eq-pencilscounting}) corresponding to 
the case $N=\{3, 9, 15, 21, 27\}$, $M=\{4, 10, 16, 22, 28\}$
are given in Tables \ref{patterns-13-2-1}, \ref{patterns-13-2-0}, 
\ref{table-13-1-2}, \ref{table-13-2-2}, respectively.
We will refer to them as the solutions $\#1$ and $\#2$ of the second group.

\begin{lemma}\label{lemma-d}
$\mathcal{D}$ is the only Cameron-Liebler line class corresponding to the solution $\#2$ 
of the second group.
\end{lemma}
\proof Let $\mathcal{L}$ be a Cameron-Liebler line class with parameter $x=13$ in $\PG(3,5)$ 
corresponding to the solution $\#2$ of the second group.
It follows from Tables \ref{patterns-13-2-1}, \ref{patterns-13-2-0}, 
\ref{table-13-1-2} that the set $\mathcal{O}$ of points of weight $3$
forms a $26$-cap, i.e., an elliptic quadric in $\PG(3,5)$ (see \cite{Barlotti,Panella}), 
while the ${26 \choose 2}=325$ lines with pattern $T_8$ are the secant lines to $\mathcal{O}$,
and the $156=26\cdot (3+3)$ lines with pattern $T_7$ or $T_1$ are the tangent lines to $\mathcal{O}$, 
moreover, half of the tangents at each point of $\mathcal{O}$ belongs to $\mathcal{L}$.
By \cite[Theorem~5.1,~Theorem~5.9]{DrudgeThesis}, $\overline{\mathcal{L}}$ is projectively equivalent 
to the line class $\mathcal{D}$, which was constructed by Bruen and Drudge in 
\cite{BruenDrudge}, \cite[Section~5]{DrudgeThesis}.
\wbull

\begin{lemma}\label{lemma-p}
$\mathcal{P}$ is the only Cameron-Liebler line class corresponding to the solution $\#1$ 
of the second group.
\end{lemma}
\proof Let $\mathcal{L}$ be a Cameron-Liebler line class with parameter $x=13$ in $\PG(3,5)$ 
corresponding to the solution $\#1$ of the second group. 
It follows from Tables \ref{patterns-13-2-1}, \ref{patterns-13-2-0}, 
\ref{table-13-1-2} that the set of planes of weight $4$ 
plus the only plane $\pi$ of weight $28$ forms a dual $26$-cap $\mathcal{O}'$, 
since every line of $\PG(3,5)$ is contained in at most two planes from $\mathcal{O}'$.
More precisely, the $300$ lines with pattern $T_{10}$ or $T_{11}$ and 
the $25$ lines with pattern $T_3$ are the secant lines to $\mathcal{O}'$, 
and the $156$ lines with pattern $T_1$, $T_4$, $T_7$ or $T_{12}$
are the tangent lines to $\mathcal{O}'$.
Now the plane $\pi$ has the following property: there exists a point $P\in \pi$ 
such that $\Li(\pi)\setminus \St(P)\subset \mathcal{L}$ and 
$(\St(P)\setminus \Li(\pi))\cap \mathcal{L}=\emptyset$.
By \cite[Lemma~1]{GMP}, the switched line class 
$\mathcal{L}^{\star}:=\mathcal{L}\cup (\St(P) \setminus \Li(\pi)) \setminus (\Li(\pi) \setminus \St(P))$
is a Cameron-Liebler line class with the same parameter $x=13$.
Moreover, $\PG(3,5)$ contains precisely $26$ planes of weight $3$ 
with respect to $\mathcal{L}^{\star}$, and hence, by Lemma \ref{lemma-d}, 
$\mathcal{L}^{\star}$ is projectively equivalent to a line class which is dual to $\mathcal{D}$.
\wbull

Theorem \ref{theo-main} follows from Lemmas \ref{lemma-group1}, \ref{lemma-12}, 
\ref{lemma-d}, \ref{lemma-p} and \cite[Theorem~5.1]{GavrilyukMetsch}.

\section{Concluding remarks}\label{section-comments}

In this section, we present some results 
regarding Cameron-Liebler line classes with parameter $x$ in $\PG(3,q)$ 
for $x$ satisfying $q+1\leq x\leq \lfloor \frac{q^2+1}{2}\rfloor$ 
(see \cite{Metsch,Metsch2}) and $q\in \{7,8\}$.

In $\PG(3,7)$, Eq. (\ref{eqn_main}) has a solution if 
$x\in \{8,9,10,13,16,17,18,21,24,25\}$.
For $x\in \{8,9,10,13\}$, there are no admissible patterns. 
For $x\in \{16,17\}$, the system of Eqs. (\ref{eq-sum-n})-(\ref{eq-pencilscounting}) 
has no solution in non-negative integers.
For $x=18$, the system of Eqs. (\ref{eq-sum-n})-(\ref{eq-pencilscounting}) 
has only 5 solutions in non-negative integers 
(the search implemented with the aid of the MIP solver Gurobi 
required 16319 seconds of CPU time on a single core of Intel Core i5-3360M 2.8 GHz processor).
We are not aware of the existence of Cameron-Liebler line classes 
with parameter 18 in $\PG(3,7)$, and so these solutions represent the smallest open case.
The solutions can be found in \cite{Matkinwww}. 
For $x\in \{21,24,25\}$, the system of Eqs. (\ref{eq-sum-n})-(\ref{eq-pencilscounting}) 
has too many unknowns relative to its rank (138 vs. 45 for $x=21$, 
178 vs. 58 for $x=24$, 190 vs. 62 for $x=25$).

In $\PG(3,8)$, Eq. (\ref{eqn_main}) has a solution if 
$x\in \{9,10,11,13,16,18,19,20,22,25,27,28,29,31\}$.
For $x\in \{9,10,13,16\}$, there are no admissible patterns. 
For $x\in \{11,18,19,20\}$, the system of Eqs. (\ref{eq-sum-n})-(\ref{eq-pencilscounting}) 
has no solution in non-negative integers 
(for $x=20$ the search implemented with the aid of the MIP solver Gurobi 
required 175442 seconds of CPU time).
For $x\in \{22,25,27,28,29,31\}$, the system of Eqs. (\ref{eq-sum-n})-(\ref{eq-pencilscounting}) 
has too many unknowns relative to its rank.

Finally, we would like to close our paper with the following open problem.
In $\PG(3,9)$, there may exist particularly interesting examples of Cameron-Liebler 
line classes related to the so-called $(a,b)$-sets of points in $\PG(n,q)$. 
We say that a set $\mathcal{O}$ of points in $\PG(n,q)$ 
such that neither $\mathcal{O}$ nor its complement is empty, 
a point or a hyperplane, is an $(a,b)$-set, $0<a\leq b<q+1$, 
if every line contains exactly $a$ or $b$ points of $\mathcal{O}$.
No examples of such sets $\mathcal{O}$ are known if $n\geq 3$.
It was shown by Tallini-Scafati, see \cite{TS,CK}, 
that if $\mathcal{O}$ exists in $\PG(n,q)$, $n\geq 3$, then 
$q$ is odd and a square and 
\begin{align*}
a &=\frac{1}{2}(q+1-\sqrt{q}(1-\epsilon)),\\
b &=\frac{1}{2}(q+1+\sqrt{q}(1+\epsilon)),\\
|\mathcal{O}| &=\frac{1}{2}(1+\frac{q^{n-1}-1}{q-1}(q+\epsilon\sqrt{q})+\delta\sqrt{q^{n-1}}),
\end{align*}
where $\epsilon=\pm 1$, $\delta=\pm 1$ 
(the complement to $\mathcal{O}$ corresponds to replacing $\epsilon$ by $-\epsilon$ and 
$\delta$ by $-\delta$). 
Thus, the first (open) case where such a set $\mathcal{O}$ may exist is $\PG(3,9)$ 
and $(a,b)=(2,5)$ (we take $\epsilon=-1$).
Further, in the case $n=3$, it follows from \cite[Lemma~8]{Penttila} that the set of lines 
that intersect $\mathcal{O}$ in the same number points is a Cameron-Liebler line class.
In $\PG(3,9)$, the corresponding line class should have parameter $32$ (if $\delta=-1$) 
or $\frac{9^2+1}{2}=41$ (if $\delta=1$), 
and one can easily determine its patterns and the numbers of lines with given patterns.
We are not aware of the existence of such line classes, 
however, all their possible planar sections are the $(2,5)$-sets in $\PG(2,9)$, 
which were found in \cite{25}. Thus, one can try to apply the approach 
from the proof of Lemma \ref{lemma-12}, which, however, would probably require 
too much computational effort.

\noindent{\bf Acknowledgments.}
The research of A.L.G. was supported by BK21plus Center for Math Research and Education at PNU.
I.M. was supported by the Russian Foundation for Basic Research (grant 17-301-50004).
Part of the work was done while I.M. was visiting the University of Science and 
Technology of China, he thanks Professor Jack Koolen for his hospitality.

\newpage
\section*{Appendix: Tables}

\begin{table}[h]
\caption{$x=13$, the first group of solutions: patterns of lines of $\mathcal{L}$}
\label{patterns-13-1-1}
\arraycolsep0.5ex
\centering
\[\small
T_1=
\begin{pmatrix}
0&  0& 0&  1& 2&  2\cr
2&  2& 2&  3& 4&  4\cr  
2&  2& 2&  3&  4&  4\cr  
2&  2& 2&  3&  4&  4\cr  
3&  3&  3&  4&  5&  5\cr  
3&  3&  3&  4&  5&  5  
\end{pmatrix}, 
T_2=
\begin{pmatrix}
0&  0&  1&  1&  1&  2\cr
1&  1&  2&  2&  2&  3\cr 
2&  2&  3&  3&  3&  4\cr  
3&  3&  4&  4&  4&  5\cr  
3&  3&  4&  4&  4&  5\cr  
3&  3&  4&  4&  4&  5 
\end{pmatrix},
T_3=
\begin{pmatrix}
0&  1&  1&  3&  3&  3\cr  
1&  2&  2&  4&  4&  4\cr  
1&  2&  2&  4&  4&  4\cr  
1&  2&  2&  4&  4&  4\cr  
1&  2&  2&  4&  4&  4\cr  
2&  3&  3&  5&  5&  5  
\end{pmatrix},
T_4=
\begin{pmatrix}
0&  1&  2&  2&  3&  3\cr  
0&  1&  2&  2&  3&  3\cr  
1&  2&  3&  3&  4&  4\cr  
1&  2&  3&  3&  4&  4\cr  
2&  3&  4&  4&  5&  5\cr  
2&  3&  4&  4&  5&  5 
\end{pmatrix},
\]
\[\small
T_5=
\begin{pmatrix}
0&  2&  2&  2&  2&  3\cr  
0&  2&  2&  2&  2&  3\cr  
0&  2&  2&  2&  2&  3\cr  
2&  4&  4&  4&  4&  5\cr  
2&  4&  4&  4&  4&  5\cr  
2&  4&  4&  4&  4&  5 
\end{pmatrix},
T_6=
\begin{pmatrix}
1&  1&  1&  2&  3&  3\cr  
1&  1&  1&  2&  3&  3\cr  
1&  1&  1&  2&  3&  3\cr  
3&  3&  3&  4&  5&  5\cr  
3&  3&  3&  4&  5&  5\cr  
3&  3&  3&  4&  5&  5 
\end{pmatrix},
T_7=
\begin{pmatrix}
1&  1&  3&  4&  4&  4\cr  
1&  1&  3&  4&  4&  4\cr  
1&  1&  3&  4&  4&  4\cr  
1&  1&  3&  4&  4&  4\cr  
1&  1&  3&  4&  4&  4\cr  
1&  1&  3&  4&  4&  4 
\end{pmatrix},
T_8=
\begin{pmatrix}
1&  2&  2&  3&  4&  5\cr  
1&  2&  2&  3&  4&  5\cr  
1&  2&  2&  3&  4&  5\cr  
1&  2&  2&  3&  4&  5\cr  
1&  2&  2&  3&  4&  5\cr  
1&  2&  2&  3&  4&  5 
\end{pmatrix}.
\]
\end{table}

\begin{table}[h]
\caption{$x=13$, the first group of solutions: patterns of lines of $\overline{\mathcal{L}}$}
\label{patterns-13-1-0}
\arraycolsep0.5ex
\centering
\[\small
T_9=
\begin{pmatrix}
0& 0& 0& 0& 0& 0\cr
1& 1& 1& 1& 1& 1\cr  
2& 2& 2& 2& 2& 2\cr  
3& 3& 3& 3& 3& 3\cr  
3& 3& 3& 3& 3& 3\cr  
4& 4& 4& 4& 4& 4 
\end{pmatrix},
T_{10}=
\begin{pmatrix}
0& 0& 0& 1& 2& 3\cr  
1& 1& 1& 2& 3& 4\cr  
1& 1& 1& 2& 3& 4\cr  
1& 1& 1& 2& 3& 4\cr  
2& 2& 2& 3& 4& 5\cr  
2& 2& 2& 3& 4& 5  
\end{pmatrix},
T_{11}=
\begin{pmatrix}
0& 0& 0& 2& 2& 2\cr  
0& 0& 0& 2& 2& 2\cr  
1& 1& 1& 3& 3& 3\cr  
2& 2& 2& 4& 4& 4\cr  
2& 2& 2& 4& 4& 4\cr  
2& 2& 2& 4& 4& 4 
\end{pmatrix},
T_{12}=
\begin{pmatrix}
0& 0& 0& 2& 2& 2\cr 
1& 1& 1& 3& 3& 3\cr 
1& 1& 1& 3& 3& 3\cr 
1& 1& 1& 3& 3& 3\cr 
1& 1& 1& 3& 3& 3\cr 
3& 3& 3& 5& 5& 5 
\end{pmatrix},
\]
\[\small
T_{13}=
\begin{pmatrix}
0& 0& 1& 1& 1& 3\cr 
0& 0& 1& 1& 1& 3\cr 
1& 1& 2& 2& 2& 4\cr 
2& 2& 3& 3& 3& 5\cr 
2& 2& 3& 3& 3& 5\cr 
2& 2& 3& 3& 3& 5
\end{pmatrix},
T_{14}=
\begin{pmatrix}
0& 0& 1& 1& 2& 2\cr 
0& 0& 1& 1& 2& 2\cr 
1& 1& 2& 2& 3& 3\cr 
1& 1& 2& 2& 3& 3\cr 
2& 2& 3& 3& 4& 4\cr 
3& 3& 4& 4& 5& 5 
\end{pmatrix},
T_{15}=
\begin{pmatrix}
0& 1& 1& 1& 1& 2\cr 
0& 1& 1& 1& 1& 2\cr 
0& 1& 1& 1& 1& 2\cr 
2& 3& 3& 3& 3& 4\cr 
2& 3& 3& 3& 3& 4\cr 
3& 4& 4& 4& 4& 5 
\end{pmatrix},
T_{16}=
\begin{pmatrix}
1& 1& 1& 1& 1& 1\cr 
1& 1& 1& 1& 1& 1\cr 
1& 1& 1& 1& 1& 1\cr 
2& 2& 2& 2& 2& 2\cr 
4& 4& 4& 4& 4& 4\cr 
4& 4& 4& 4& 4& 4
\end{pmatrix}.
\]
\end{table}

\begin{table}[]
\centering
\caption{The feasible weight distributions for the first group of solutions in the case $x=13$}
\label{table-13-1}
\begin{tabular}{|l|l|l|}
\hline
\# & $(n_0, n_6, n_{12}, n_{18}, n_{24}, n_{30})$ & $(m_1, m_7, m_{13}, m_{19}, m_{25}, m_{31})$ \\ \hline
1  &  $(0, 37, 19, 72, 28, 0)$                    & $(0, 28, 72, 19, 37, 0)$                     \\ \hline
2  &  $(0, 33, 31, 60, 32, 0)$                    & $(0, 32, 60, 31, 33, 0)$                     \\ \hline
3  &  $(1, 31, 31, 62, 31, 0)$                    & $(0, 31, 62, 31, 31, 1)$                     \\ \hline
\end{tabular}
\end{table}

\begin{table}[]
\centering
\caption{The first group of solutions of Eqs. (\ref{eq-allsum})-(\ref{eq-pencilscounting}) 
in the case $x=13$}
\label{table-13-2}
\begin{tabular}{|c|c|c|c|c|c|c|c|c|c|c|c|c|c|c|c|c|}
\hline
\multicolumn{1}{|l|}{} & \multicolumn{8}{c|}{$\mathcal{L}$}                    & \multicolumn{8}{c|}{$\overline{\mathcal{L}}$}                               \\ \hline
\# & $z_1$ & $z_2$ & $z_3$ & $z_4$ & $z_5$ & $z_6$ & $z_7$ & $z_8$ & $z_9$ & $z_{10}$ & $z_{11}$ & $z_{12}$ & $z_{13}$ & $z_{14}$ & $z_{15}$ & $z_{16}$\\ \hline
1 & 
198 & 24  &  111 & 18 & 1 & 18 & 33 & 0 & 0  & 24 & 18 & 1 & 198 & 18 & 111 & 33 \\ \hline
2 & 
112 & 86  &  85  & 96 & 1 &  2 & 21 & 0 & 0  & 86 & 2  & 1 & 112 & 96 & 85  &  21 \\ \hline
3 & 
93  & 93  &  93  & 93 & 0 &  0 & 0 & 31 & 31 & 93 & 0  & 0 & 93  & 93 & 93  &  0\\ \hline
\end{tabular}
\end{table}

\begin{table}[h]
\caption{$x=12$: patterns of lines of $\mathcal{L}'$}
\label{patterns-12-1}
\arraycolsep0.5ex
\centering
\[\small
T_1=
\begin{pmatrix}
0 & 0 & 0 & 1 & 2 & 2\cr 
1 & 1 & 1 & 2 & 3 & 3\cr 
2 & 2 & 2 & 3 & 4 & 4\cr 
2 & 2 & 2 & 3 & 4 & 4\cr 
3 & 3 & 3 & 4 & 5 & 5\cr 
3 & 3 & 3 & 4 & 5 & 5 
\end{pmatrix},
T_2=
\begin{pmatrix}
0 & 0 & 1 & 1 & 1 & 2\cr 
1 & 1 & 2 & 2 & 2 & 3\cr
1 & 1 & 2 & 2 & 2 & 3\cr 
3 & 3 & 4 & 4 & 4 & 5\cr 
3 & 3 & 4 & 4 & 4 & 5\cr 
3 & 3 & 4 & 4 & 4 & 5 
\end{pmatrix},
T_3=
\begin{pmatrix}
0 & 1 & 1 & 3 & 3 & 3\cr 
0 & 1 & 1 & 3 & 3 & 3\cr 
1 & 2 & 2 & 4 & 4 & 4\cr 
1 & 2 & 2 & 4 & 4 & 4\cr 
1 & 2 & 2 & 4 & 4 & 4\cr 
2 & 3 & 3 & 5 & 5 & 5 
\end{pmatrix},
T_4=
\begin{pmatrix}
0 & 1 & 2 & 2 & 3 & 3\cr 
0 & 1 & 2 & 2 & 3 & 3\cr 
0 & 1 & 2 & 2 & 3 & 3\cr 
1 & 2 & 3 & 3 & 4 & 4\cr 
2 & 3 & 4 & 4 & 5 & 5\cr 
2 & 3 & 4 & 4 & 5 & 5
\end{pmatrix}.
\]
\end{table}

\begin{table}[h]
\caption{$x=12$: patterns of lines of $\overline{\mathcal{L}'}$}
\label{patterns-12-0}
\arraycolsep0.5ex
\centering
\[\small
T_5=
\begin{pmatrix}
0 & 0 & 0 & 0 & 0 & 0\cr 
1 & 1 & 1 & 1 & 1 & 1\cr 
2 & 2 & 2 & 2 & 2 & 2\cr 
2 & 2 & 2 & 2 & 2 & 2\cr 
3 & 3 & 3 & 3 & 3 & 3\cr 
4 & 4 & 4 & 4 & 4 & 4 
\end{pmatrix},
T_6=
\begin{pmatrix}
0 & 0 & 0 & 1 & 2 & 3\cr 
1 & 1 & 1 & 2 & 3 & 4\cr 
1 & 1 & 1 & 2 & 3 & 4\cr 
1 & 1 & 1 & 2 & 3 & 4\cr 
1 & 1 & 1 & 2 & 3 & 4\cr
2 & 2 & 2 & 3 & 4 & 5 
\end{pmatrix},
T_7=
\begin{pmatrix}
0 & 0 & 1 & 1 & 1 & 3\cr 
0 & 0 & 1 & 1 & 1 & 3\cr 
1 & 1 & 2 & 2 & 2 & 4\cr 
1 & 1 & 2 & 2 & 2 & 4\cr 
2 & 2 & 3 & 3 & 3 & 5\cr 
2 & 2 & 3 & 3 & 3 & 5 
\end{pmatrix},
\]
\[\small
T_8=
\begin{pmatrix}
0 & 0 & 1 & 1 & 2 & 2\cr 
0 & 0 & 1 & 1 & 2 & 2\cr 
1 & 1 & 2 & 2 & 3 & 3\cr 
1 & 1 & 2 & 2 & 3 & 3\cr 
1 & 1 & 2 & 2 & 3 & 3\cr 
3 & 3 & 4 & 4 & 5 & 5 
\end{pmatrix},
T_9=
\begin{pmatrix}
0 & 1 & 1 & 1 & 1 & 2\cr 
0 & 1 & 1 & 1 & 1 & 2\cr 
0 & 1 & 1 & 1 & 1 & 2\cr 
1 & 2 & 2 & 2 & 2 & 3\cr 
2 & 3 & 3 & 3 & 3 & 4\cr 
3 & 4 & 4 & 4 & 4 & 5 
\end{pmatrix},
T_{10}=
\begin{pmatrix}
0 & 1 & 2 & 2 & 3 & 4\cr 
0 & 1 & 2 & 2 & 3 & 4\cr 
0 & 1 & 2 & 2 & 3 & 4\cr 
0 & 1 & 2 & 2 & 3 & 4\cr 
0 & 1 & 2 & 2 & 3 & 4\cr 
0 & 1 & 2 & 2 & 3 & 4 
\end{pmatrix}.
\]
\end{table}

\begin{table}[]
\centering
\caption{The solution of Eqs. (\ref{eq-allsum})-(\ref{eq-pencilscounting}) 
in the case $x=12$}
\label{table-12}
\begin{tabular}{|c|c|c|c|c|c|c|c|c|c|c|}
\hline
\multicolumn{1}{|l|}{} & \multicolumn{4}{c|}{$\mathcal{L}'$}                    & \multicolumn{6}{c|}{$\overline{\mathcal{L}'}$}                               \\ \hline
\# & $z_1$ & $z_2$ & $z_3$ & $z_4$ & $z_5$ & $z_6$ & $z_7$ & $z_8$ & $z_9$ & $z_{10}$  
\\ \hline
1 & 
93 & 93 & 93 & 93 & 31 & 93 & 93 & 93 & 93 & 31 \\ \hline
\end{tabular}
\end{table}

\begin{table}[h]
\caption{$x=13$, the second group of solutions: patterns of lines of $\mathcal{L}$}
\label{patterns-13-2-1}
\arraycolsep0.5ex
\centering
\[\small
T_1=
\begin{pmatrix}
0 & 0 & 0 & 0 & 0 & 2\cr 
3 & 3 & 3 & 3 & 3 & 5\cr  
3 & 3 & 3 & 3 & 3 & 5\cr 
3 & 3 & 3 & 3 & 3 & 5\cr
3 & 3 & 3 & 3 & 3 & 5\cr
3 & 3 & 3 & 3 & 3 & 5
\end{pmatrix}, 
T_2=
\begin{pmatrix}
0 & 1 & 1 & 2 & 2 & 2\cr 
1 & 2 & 2 & 3 & 3 & 3\cr
1 & 2 & 2 & 3 & 3 & 3\cr
1 & 2 & 2 & 3 & 3 & 3\cr
3 & 4 & 4 & 5 & 5 & 5\cr
3 & 4 & 4 & 5 & 5 & 5
\end{pmatrix},
T_3=
\begin{pmatrix}
0 & 2 & 2 & 3 & 3 & 4\cr
0 & 2 & 2 & 3 & 3 & 4\cr
0 & 2 & 2 & 3 & 3 & 4\cr
1 & 3 & 3 & 4 & 4 & 5\cr
1 & 3 & 3 & 4 & 4 & 5\cr
1 & 3 & 3 & 4 & 4 & 5
\end{pmatrix},
\]
\[\small
T_4=
\begin{pmatrix}
0 & 2 & 3 & 3 & 3 & 3\cr
0 & 2 & 3 & 3 & 3 & 3\cr
0 & 2 & 3 & 3 & 3 & 3\cr
0 & 2 & 3 & 3 & 3 & 3\cr
1 & 3 & 4 & 4 & 4 & 4\cr 
2 & 4 & 5 & 5 & 5 & 5
\end{pmatrix},
T_5=
\begin{pmatrix}
1 & 1 & 1 & 1 & 2 & 2\cr
1 & 1 & 1 & 1 & 2 & 2\cr
2 & 2 & 2 & 2 & 3 & 3\cr
3 & 3 & 3 & 3 & 4 & 4\cr
4 & 4 & 4 & 4 & 5 & 5\cr
4 & 4 & 4 & 4 & 5 & 5
\end{pmatrix},
T_6=
\begin{pmatrix}
1 & 1 & 2 & 2 & 4 & 4\cr
1 & 1 & 2 & 2 & 4 & 4\cr
1 & 1 & 2 & 2 & 4 & 4\cr
2 & 2 & 3 & 3 & 5 & 5\cr
2 & 2 & 3 & 3 & 5 & 5\cr
2 & 2 & 3 & 3 & 5 & 5
\end{pmatrix}.
\]
\end{table}

\begin{table}[h]
\caption{$x=13$, the second group of solutions: patterns of lines of $\overline{\mathcal{L}}$}
\label{patterns-13-2-0}
\arraycolsep0.5ex
\centering
\[\small
T_7=
\begin{pmatrix}
0 & 0 & 0 & 0 & 0 & 3\cr
2 & 2 & 2 & 2 & 2 & 5\cr
2 & 2 & 2 & 2 & 2 & 5\cr
2 & 2 & 2 & 2 & 2 & 5\cr
2 & 2 & 2 & 2 & 2 & 5\cr
2 & 2 & 2 & 2 & 2 & 5
\end{pmatrix},
T_{8}=
\begin{pmatrix}
0 & 0 & 0 & 1 & 1 & 1\cr
0 & 0 & 0 & 1 & 1 & 1\cr
2 & 2 & 2 & 3 & 3 & 3\cr
2 & 2 & 2 & 3 & 3 & 3\cr
3 & 3 & 3 & 4 & 4 & 4\cr
3 & 3 & 3 & 4 & 4 & 4
\end{pmatrix},
T_{9}=
\begin{pmatrix}
0 & 0 & 0 & 1 & 1 & 1\cr
1 & 1 & 1 & 2 & 2 & 2\cr
1 & 1 & 1 & 2 & 2 & 2\cr
2 & 2 & 2 & 3 & 3 & 3\cr
2 & 2 & 2 & 3 & 3 & 3\cr
4 & 4 & 4 & 5 & 5 & 5
\end{pmatrix},
\]
\[\small
T_{10}=
\begin{pmatrix}
0 & 0 & 1 & 2 & 3 & 3\cr
0 & 0 & 1 & 2 & 3 & 3\cr
1 & 1 & 2 & 3 & 4 & 4\cr
1 & 1 & 2 & 3 & 4 & 4\cr
1 & 1 & 2 & 3 & 4 & 4\cr
1 & 1 & 2 & 3 & 4 & 4
\end{pmatrix},
T_{11}=
\begin{pmatrix}
0 & 0 & 2 & 2 & 2 & 3\cr
0 & 0 & 2 & 2 & 2 & 3\cr
0 & 0 & 2 & 2 & 2 & 3\cr
1 & 1 & 3 & 3 & 3 & 4\cr
1 & 1 & 3 & 3 & 3 & 4\cr
2 & 2 & 4 & 4 & 4 & 5
\end{pmatrix},
T_{12}=
\begin{pmatrix}
0 & 1 & 2 & 2 & 2 & 2\cr
0 & 1 & 2 & 2 & 2 & 2\cr
0 & 1 & 2 & 2 & 2 & 2\cr
0 & 1 & 2 & 2 & 2 & 2\cr
1 & 2 & 3 & 3 & 3 & 3\cr
3 & 4 & 5 & 5 & 5 & 5
\end{pmatrix}.
\]
\end{table}

\begin{table}[]
\centering
\caption{The feasible weight distributions for the second group of solutions in the case $x=13$}
\label{table-13-1-2}
\begin{tabular}{|l|l|l|}
\hline
\# & $(n_3, n_9, n_{15}, n_{21}, n_{27})$ & $(m_4, m_{10}, m_{16}, m_{22}, m_{28})$ \\ \hline
1  &  $(1,  50,  65,  15,  25)$  & $(25, 15, 65, 50, 1)$                     \\ \hline
2  &  $(26, 0,  65,  65,  0)$   & $(0, 65, 65,  0, 26)$                     \\ \hline
\end{tabular}
\end{table}

\begin{table}[]
\centering
\caption{The second group of solutions of Eqs. (\ref{eq-allsum})-(\ref{eq-pencilscounting}) 
in the case $x=13$}
\label{table-13-2-2}
\begin{tabular}{|c|c|c|c|c|c|c|c|c|c|c|c|c|}
\hline
\multicolumn{1}{|l|}{} & \multicolumn{6}{c|}{$\mathcal{L}$}                    & \multicolumn{6}{c|}{$\overline{\mathcal{L}}$}                               \\ \hline
\# & $z_1$ & $z_2$ & $z_3$ & $z_4$ & $z_5$ & $z_6$ & $z_7$ & $z_8$ & $z_9$ & $z_{10}$ & $z_{11}$ & $z_{12}$ 
\\ \hline
1 & 
3 & 150 & 25 & 75 & 150 & 0 & 3 & 0 & 25 & 150 & 150 & 75 \\ \hline
2 & 
78 & 0  &  0  & 0 & 0 &  325 & 78 & 325 & 0  & 0 & 0  & 0 \\ \hline
\end{tabular}
\end{table}

\end{document}